\documentclass[referee,sn-basic]{sn-jnl}
\usepackage{graphicx}
\usepackage{graphicx} 
\usepackage{amsmath}
\usepackage{amssymb}
\usepackage{amsfonts}
\usepackage[title]{appendix}
\usepackage{xcolor}
\usepackage{textcomp}%
\usepackage{manyfoot}%
\usepackage{booktabs}%
\usepackage{algorithm}%
\usepackage{algorithmicx}%
\usepackage{algpseudocode}%
\usepackage{listings}%
\usepackage{physics}
\usepackage{float}
\usepackage{multirow}
\usepackage{epstopdf}
\usepackage{subcaption}
\begin{document}

\title{SIR Epidemics in Interconnected Networks: threshold curve and phase transition}

\author*[1]{\fnm{Saswata} \sur{Das}}\email{ersaswatadas@ksu.edu}
\author*[1]{\fnm{Mohammad Hossein} \sur{Samaei}}\email{msamaei@ksu.edu}
\author[1]{\fnm{Caterina} \sur{Scoglio}}\email{caterina@ksu.edu}

\affil[1]{\orgdiv{Electrical and Computer Engineering}, \orgname{Kansas State University}, \orgaddress{ \city{Manhattan}, \state{Kansas}, \postcode{66502}, \country{United States of America}}}

\abstract
{For simplicity of mathematical modeling of epidemic spreading, the assumption is that hosts have identical rates of disease-causing contacts. However, in the real world, the scenario is different. The network-based framework allows us to capture the complex interdependencies and structural heterogeneity present in real-world systems. 
We examine two distinct scenarios involving the dynamics of Susceptible-Infected-Recovered (SIR) in interconnected networks. In the first part, we show how the epidemic threshold of a contact network changes as a result of being coupled with another network for a fixed infection strength. The model employed in this work considers both the contact networks and interconnections as generic. We have depicted the epidemic threshold curve for interconnected networks, considering the assumption that the infection could be initially present in either one or both of the networks. If the normalized infection strengths are above the threshold curve, the infection spreads, whereas if the normalized infection strengths are below the threshold curve, the disease does not spread. This is true for any level of interconnection. 
In the second part, we investigate the spillover phenomenon, where the disease in a novel host population network comes from a reservoir network. We have observed a clear phase transition when the number of links or the inter-network infection rate exceeds a certain threshold, keeping all other parameters constant. We observe two regimes for spillover: a major spillover region and a minor spillover region based on interpopulation links (fraction of links between two networks) and inter-network infection strength (infection rate between reservoir and host network). If the interpopulation links and inter-network infection strength are in the major spillover region, the spillover probability is high, while if the former parameters are in the minor spillover region, the spillover probability is low. When the number of infected individuals within a reservoir network is nearly equal, and the inter-network infection strength remains constant, the threshold number of links required to achieve the spillover threshold condition varies based on the network topology. 
Overall, this work contributes to the understanding of SIR dynamics in interconnected networks and sheds light on the behavior of epidemics in complex systems.} 


\maketitle

\section{Introduction}\label{sec1}
Various types of contagious diseases are found all over the world. It is crucial to prevent the spread of these diseases. Consequently, extensive research is being conducted to halt their transmission. In addition to other research endeavors, scientists are also crafting diverse mathematical models. These models can aid in understanding disease-spreading mechanisms, estimating potential outbreak sizes, and formulating more effective mitigation strategies.

Considerable progress has been made in the realm of mathematical research on epidemic diseases, encompassing both theoretical and practical applications. Numerous epidemic models are formulated as dynamical systems of ordinary differential equations. Delay differential equations have demonstrated their efficacy in capturing the variability observed in infectious periods within diverse epidemic models. Furthermore, incorporating factors like age structure and spatial considerations has contributed to the development of partial differential equations \citep{ji2014threshold}.

Compartmental models are also very useful for studying epidemic diseases. These models for epidemics are of different types, such as the SIS (Susceptible-Infected-Susceptible) model, the SIR (Susceptible-Infected-Recovered) model, and the SEIR (Susceptible-Exposed-Infected-Recovered) model, among others \citep{brauer2008compartmental}. In the SIR model, infected agents recover and become permanently immune, never becoming infected again. The SIR model is particularly effective in explaining diseases like influenza, COVID-19, and other contagious diseases. 

Many epidemiological models simplify the assumptions regarding the patterns of disease-causing interactions among hosts. Specifically, in homogeneous-mixing models, it is assumed that hosts have identical rates of disease-causing contact. In recent years, several network-based approaches have been developed to explicitly model the heterogeneity in host contact patterns \citep {bansal2007individual}.

A contact network plays a crucial role in facilitating the transmission of epidemics. This network is composed of nodes, which correspond to individual agents, and links, which signify the connections between any two individuals. At its simplest, a contact network takes on a binary form with two distinct values: a non-zero value denoting contact between a pair of individuals and 0 indicating no contact. Each node within this network represents an individual, while the interconnecting links between nodes provide insights into the quantity and types of interactions. For instance, these links shed light on an individual's associations with susceptible or infected individuals, thereby offering a comprehensive view of the underlying dynamics of disease transmission within the network.

In many cases, a single generic network is used to model epidemic diseases, but in reality, every network has some interaction with other networks. That's why the concept of interconnected networks is a very critical research area nowadays, especially for epidemic disease modeling \citep{wang2014epidemic,wang2011effects}.  

It is a proven result that in a generic contact network, if the ratio of infection to cure rate is less than a particular threshold, the infection cannot survive in the network. However, for SIS spreading dynamics in two interconnected networks, the epidemic can survive though the infection-to-cure rate in both networks is less than the epidemic threshold \citep{sahneh2013effect}. In that case, we have seen an epidemic threshold curve, which is plotted to show the epidemic threshold of one network as a function of effective infection strength (infection to cure rate) in the other network for a specific level of interconnection. To our best knowledge, this result has not been proven for SIR spreading dynamics.     

In the initial portion of this paper, we have demonstrated the alteration of the epidemic threshold of a specific network due to its coupling with another network, assuming a constant infection strength.  However, we have specifically examined the change in the epidemic threshold concerning the dynamics of SIR spreading. We have plotted three distinct epidemic threshold curves corresponding to three different types of coupling: weak, medium, and strong. We show that for SIR dynamics, although the epidemic threshold of both networks is less than the epidemic threshold, the infection can still survive in the whole interconnected system.

In our work in depicting the epidemic threshold curve for interconnected networks, the consideration has been that few initial infections either die out or invade both networks, causing an epidemic.

However, in reality, there are certain diseases where the infection is endemically present in one population represented by a contact network, and from that, the disease spills to another population represented by another contact network if the two networks are interconnected. This event is known as spillover.

In the second portion, we have done significant work on spillover in various network-based models and obtained a clear phase transition for the probability of spillover. These network-based models for spillover include both homogeneous and heterogeneous node degree distributions. We have also drawn a spillover threshold curve based on inter-population link density and inter-network infection strength that divides the entire area into two regions: major and minor spillover. If the inter-network infection strength and interpopulation link density are in the major spillover region, the probability of spillover is very high, whereas if the inter-network infection strength and interpopulation link density are in the minor spillover region, the probability of spillover is low.

Summarizing the novel contributions of this paper are:
\begin{itemize}
  \item The determination of an epidemic threshold curve for SIR spreading in interconnected networks.
  \item The discovery of a phase transition for spillover as a function of the level of interconnection.
  \item The determination of a major or minor spillover probability region as a function of the level of interconnection and inter-network infection strength.
\end{itemize}

The rest of the paper is organized as follows. Section 2 presents some background on spillover and interconnected networks. The epidemic threshold curve for SIR dynamics in interconnected networks is illustrated in Section 3. Modeling spillover and discovering the phase transition through extensive simulations are shown in Section 4.

\section{Background}\label{sec2}

Multi-layer networks are formed by several networks that evolve and interact with each other. These networks are ubiquitous and include social networks, financial markets, and multi-modal transportation systems. The multi-layer structure of these networks strongly affects the properties of dynamical and stochastic processes defined on them, which can display unexpected characteristics \citep{bianconi2018multilayer}.

As the field of network science has dramatically increased, researchers have found that most complex systems do not work in isolation. Every complex system depends on another complex system to some extent. For example, communication systems often depend on the power grid to be operational \citep{gao2022introduction,d2014networks}. So, if there is a failure in any part of an interconnected system, that fault propagates to other parts of the system. So, for the proper functioning of an interconnected system, insight into fault propagation mechanisms in an interconnected system is very important. This kind of research also provides insight into the robustness of the system \citep{shekhtman2023spreading}. Similarly, research has been done to revive a failed network through microscopic interventions. The complex system represented by a complex network generally fails due to node, link removal, or due to reduction in link weights. But just reversing the topological damage, i.e., reconstruction of links, nodes, or increase of link weight, does not guarantee the spontaneous recovery of a system. So scientists have come out with a two-step intervention process by which a system can be steered towards its functionality after reversing the topological damage \citep{sanhedrai2022reviving}. Radicchi et al. have shown that depending on the relative importance of inter and intra-layer connection, the entire interconnected system can have two regimes; in one regime, various layers act as independent entities, and in another regime, the whole system behaves as a single network \citep{radicchi2013abrupt}. These are some examples of various research conducted on interconnected networks.

Interconnected networks are also used to model epidemic spreading. Researchers have shown how the epidemic threshold of an interconnected network structure varies for SIS spreading dynamics without any approximation. They have also shown the upper and lower bound of the epidemic threshold and how it is related to properties of network parameters (eigenvalue and eigenvector) \citep{wang2013effect}. Scientists have found that for an interconnected network system, there exists a global threshold above which the infection prevails in every network of the subsystem and below which the infection dies out. Another finding indicated that having a diverse structure enhances the likelihood of infection, with this impact being particularly noticeable in interconnected network systems \citep{zhu2015mean}. The epidemic threshold in two interconnected networks is always lower than any of the two-component networks. Moreover, in interconnected networks, interconnection correlation has no significant contribution to epidemic size \citep{wang2012epidemics}. Dickison and colleagues discovered that when considering SIR spreading dynamics within a connected network, there are two clear modes. In the case of strongly coupled interconnected networks, there's a threshold infection strength beyond which the epidemic invades the entire interconnected structure, while below this threshold, the epidemics die away. But for weakly coupled networks, a mixed phase exists, where the epidemic does not spread in the whole interconnected system, and interconnections affect the less interconnected network mostly \citep{dickison2012epidemics}. The directed interconnected network is also used for epidemic modeling. For various types of directed interconnected networks, mathematical expression of $R_{0}$. The disease will become endemic if $R_{0}$ is greater than 1; otherwise, it will die out. Coupling can increase $R_{0}$ and even make a disease endemic. There are certain considerations for epidemic prevalence in a single sub-network, which is only possible in directed network \citep{jia2019epidemic}. A new type of model has been made to investigate the epidemics on interconnected networks. One contact network has a fixed infection rate, and another has a periodic infection rate. They have some novel findings regarding $R_{0}$ based on their model and showed the dependence of infection rate and other network parameters on $R_{0}$ \citep{xu2019propagation}.  

In recent years, many works have been done on modeling spillover. Nandi et al. have considered the effect of seasonal variation on transmission and recovery rates in the context of spillover. They have focused on the direct transmission of pathogens between humans and animals and considered all the infection and recovery rates are periodic. A branching process approximation has been applied near the disease-free equilibrium to predict the first spillover event. It also shows how the probability of spillover depends on the human-to-human infection rate, human-to-animal infection rate, and animal recovery rate \citep{nandi2021probability}. Grange et al. have identified the host, viral, and environmental risk factors contributing to zoonotic virus spillover and spread in humans. They have also developed an interactive web tool that estimates the risk score of wildlife-origin viruses and lists a number of viruses based on their risk score \citep{grange2021ranking}. Ellwanger et al. have reviewed the basic aspects and the main factors involved in zoonotic spillover. The focus was on the role of inter-species interaction, phylogenetic distance between host and species, environmental drivers, and specific characteristics of the pathogen. They have also shed light on preventing zoonotic spillover in various ways \citep{ellwanger2021zoonotic}. Royce et al. have used a new mathematical spillover model with an intermediate host. The agents are in three categories: wild animals, domestic animals, and humans. They have assumed that the pathogen in domestic animals mutates itself and becomes strong enough to infect a human. Even though $R_{0}$ is less than 1 in humans, it can still infect the human population \citep{royce2020mathematically}. A synthetic framework for animal-to-human virus transmission is proposed by Plowright et al., and this study integrates the relevant spillover mechanism. According to the authors, all zoonotic pathogens must overcome a hierarchical series of barriers to cause spillover. If any barriers are impenetrable, then the spillover cannot happen. They describe these barriers in detail and claim that the probability of spillover is determined by the interaction among the barriers and the associated bottlenecks that might prevent cross-species transmission \citep{plowright2017pathways}. Rees et al. review the mathematical models of spillover. There are two criteria for selecting diseases included in models. The first one is the disease must be zoonotic, and the second one is the pathogen must be alive for 48 hours and should be able to infect humans. This appears to describe the scope of future research in zoonotic spillover, model validation, and other important information \citep{rees2021transmission}.  We have been discussing only zoonotic spillover till now, but recently, spillover from plants has also been spotted. A 61-year-old person with no previous history of the disease has been infected by a plant fungus in Kolkata (a city in Eastern India). The scientific name of the plant pathogen is $Chondrostereum$ $purpureum$ \citep{dutta2023paratracheal}.

\section{Modeling SIR spreading in interconnected network}\label{sec3}

 One of the most commonly employed compartmental models in epidemic modeling is the SIR model, which stands for Susceptible-Infected-Recovered. The SIR model is highly popular for understanding and predicting the spread of infectious diseases within a population.

S (Susceptible): A susceptible individual is someone who has the potential to become infected when exposed to an individual carrying the infection.

I (Infected): The infected individuals represent the people who have already contracted the disease and can transmit it to susceptible individuals.

R (Recovered): Recovered individuals are those who have previously been infected with the disease and have subsequently overcome it, resulting in immunity and no longer being capable of transmitting the infection.

\begin{figure}[h]
        \centering
        \includegraphics[width=1\textwidth]{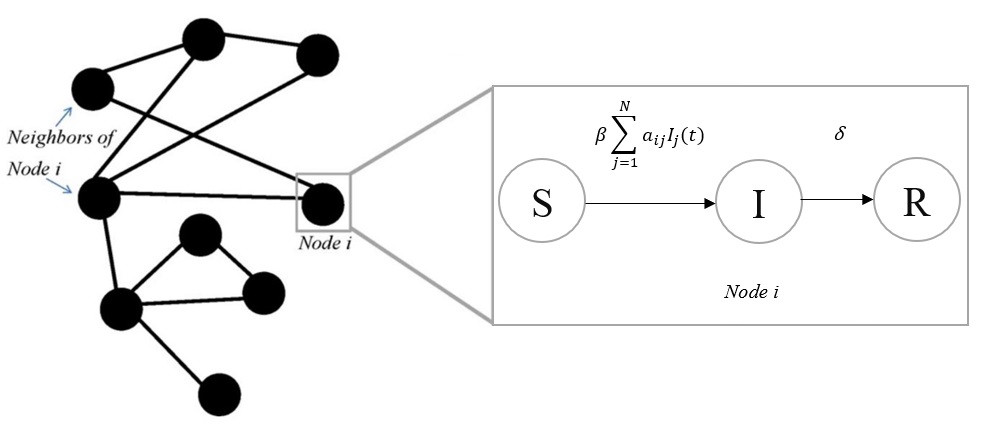}
        \caption{Schematic of a contact network along with the agent-level stochastic transition diagram for agent $i$ according to the SIR epidemic spreading model. $\beta$ and $\delta$ denote the infection and recovery rate, respectively. The notation $a_{ij}$ is used to represent the connection between agent $i$ and agent $j$. Specifically, if the value of $a_{ij}$ is 0, it signifies the absence of a link between the two agents. Conversely, if $a_{ij}$ is nonzero, it indicates the presence of a link between agent $i$ and agent $j$. $I_j(t)$ is the state probability of being infected of node $j$ at time $t$.}
        \label{fig:1}
    \end{figure}

The SIR model on a contact network with $N$ number of agents is shown in Figure \ref{fig:1}.

The equation of the SIR model with a mean-field approximation for $i^{th}$ agent in a network having $N$ number of nodes (without considering demography) is 
     \begin{align*}
        \frac{{dS_{i}(t)}}{{dt}} &= -\beta \cdot S_{i}(t) (\sum_{j=1}^{N}a_{ij}I_{j}(t)) \\
        \frac{{dI_{i}(t)}}{{dt}} &= \beta \cdot S_{i}(t) (\sum_{j=1}^{N}a_{ij}I_{j}(t)) - \delta\cdot I_{i}(t), \\
        \frac{{dR_{i}(t)}}{{dt}} &= \delta \cdot I_{i}(t).
    \end{align*}

Considering a node denoted by $i$ within a network comprising $N$ nodes, where each node is assigned labels ranging from 1 to $N$. The parameter $\beta$ is the infection rate, and $\delta$ is the recovery rate.  $S_{i}(t)$, $I_{i}(t)$, $R_{i}(t)$ denote the state probability of agent $i$ being in the susceptible, infected, or recovered compartment, respectively.  $a_{ij}$ represents an entry in the network's adjacency matrix, where the values can only be 0 or 1. A nonzero value indicates a connection between agent $i$ and agent $j$, while a value of 0 signifies the absence of a link between agent $i$ and $j$. It can also be written that at any point in time, $S_{i}(t)+I_{i}(t)+R_{i}(t)$=1.

\subsection{Generalized epidemic threshold curve}

Let us consider two groups of agents of sizes: $N_1$ and $N_2$. Let's take the first graph, $G_1$, whose agents are labeled as 1 to $N_1$, and in the second graph, the agents are labeled as $(N_1+1)$ to $(N_1+N_2)$. The adjacency matrix of graph $G_1$ is denoted by $A_{11}$, and the adjacency matrix of graph $G_2$ is denoted by $A_{22}$. The elements of sub-matrix $A_{12}$ and $A_{21}$ denotes the connection between node $i$ and node $j$ , where node $i$ belongs to graph $G_1$ and node $j$ belongs to graph $G_2$. So, the adjacency matrix for the two interconnected networks is: 
\[\begin{bmatrix}
       A_{11} & A_{12} \\
       A_{21} & A_{22} \\
\end{bmatrix}\]

The topology of networks $G_1$ and $G_2$ is undirected. As the interconnection between both networks $G_1$ and $G_2$ is undirected, it can be written as $A_{12}=(A_{21})^T$. The infection rates $\beta_{11}, \beta_{12}, \beta_{21}, \beta_{22}$ are such that a susceptible agent of graph $G_m$ receives the infection from an infected agent of graph $G_n$ with the infection rate, $\beta_{mn}$ for $m,n\in[1,2]$. The recovery rates of the agents of both networks are the same and equal to $\mu$. The equations for the change in the state probability of infection of node $i$ at time $t$ for SIR spreading dynamics in the interconnected network can be expressed in the following manner:

\begin{equation}
 \frac{dI_i(t)}{dt}= S_i(t)*[{\beta_{11}}\sum_{j=1}^{N_1} a_{ij}I_j(t)+\beta_{12}\sum_{j=N_1+1}^{N_1+N_2} a_{ij}I_j(t)]-\mu I_i(t); \; \;
 i=1,.., N_1
\end{equation}
\begin{equation}
 \frac{dI_i(t)}{dt}= S_i(t)*[{\beta_{21}}\sum_{j=1}^{N_1} a_{ij}I_j(t)+\beta_{22}\sum_{j=N_1+1}^{N_1+N_2} a_{ij}I_j(t)]-\mu I_i(t); \; \;
 i=N_1+1,.., N_1+N_2
\end{equation}

Our assumption is that if there is no interconnection, then infection cannot survive in $G_2$; that is, $\beta_{22}/\mu$ will be less than the epidemic threshold of $G_2$. So, in our case, the effective infection strength $\beta_{22}/\mu$ will be between zero to $1/\lambda(A_{22})$. Here $\lambda(A_{22})$ is the highest eigenvalue of $A_{22}$. 

As the initial fraction of infected individuals is very small, we can assume that $S_i(0)\approx 1$. So, equations (1) and (2) can be written as follows: 

\begin{equation}
 \frac{dI_i(t)}{dt}\approx[{\beta_{11}}\sum_{j=1}^{N_1} a_{ij}I_j(t)+\beta_{12}\sum_{j=N_1+1}^{N_1+N_2} a_{ij}I_j(t)]-\mu I_i(t); \; \;
 i=1,.., N_1
\end{equation}
\begin{equation}
 \frac{dI_i(t)}{dt}\approx[{\beta_{21}}\sum_{j=1}^{N_1} a_{ij}I_j(t)+\beta_{22}\sum_{j=N_1+1}^{N_1+N_2} a_{ij}I_j(t)]-\mu I_i(t); \; \;
 i=N_1+1,.., N_1+N_2
\end{equation}
To find the epidemic threshold, we follow the approach of Youssef et al. \citep{youssef2011individual}.

Equations (3) and (4) can be written as follows:
\begin{equation}
\frac{dI_i(t)}{dt}\approx\sum_{j}{L}_{ij}I_j(t).
\end{equation}

The element ${L}_{ij}$ = ($\beta a_{ij}-\mu\delta_{ij}$) of equation (5) is the element of the Jacobian matrix ${L}$ defined as follows (equation (6)).

\begin{equation}
{L} =
\begin{bmatrix}
       \beta_{11}A_{11} & \beta_{12}A_{12} \\
       \beta_{21}A_{21} & \beta_{22}A_{22} \\
\end{bmatrix}-\mu*I_{(N1+N2)(N1+N2)}
\end{equation}

where $\delta_{ij}$ is the Kronecker delta function and $I_{(N1+N2)(N1+N2)}$ is the identity matrix of dimension $(N1+N2)(N1+N2)$. We now perform the spectral analysis of the Jacobian matrix to study the growth or decay of epidemics in interconnected systems. In the above-mentioned Jacobian matrix, the elements were in the form of $\beta a_{ij}-\mu\delta_{ij}$. However, instead of using these elements in the form of $\beta a_{ij}-\mu\delta_{ij}$, we can write them as $\mu[\tau a_{ij}-\delta_{ij}]$, where $\tau=\beta/\mu$. Now a new matrix $\Bar{L}$ can be constructed where each element will be in the form $[\tau a_{ij}-\delta_{ij}]$. 

For interconnected graphs, this new matrix ($\Bar{L}$) is shown in equation (7)
\begin{equation}
\bar{L}=
\begin{bmatrix}
       \tau_{11}A_{11} & \tau_{12}A_{12} \\
       \tau_{21}A_{21} & \tau_{22}A_{22} \\
\end{bmatrix}-I_{(N1+N2)(N1+N2)}
\end{equation}
where $\tau_{11}=(\beta_{11}/\mu); \tau_{12}=(\beta_{12}/\mu) ; \tau_{21}=(\beta_{21}/\mu) ; \tau_{22}=(\beta_{22}/\mu).$
Thus, we can define \(L = \mu \bar{L}\). Let us denote the matrix \(\bar{L}\) as \([\widetilde{L} - I]\). It is known that a system with negative eigenvalues indicates that perturbations will decay over time, leading the system towards stability. Since the recovery rate, \(\mu\), is a positive quantity, it follows that if the eigenvalues of \([\widetilde{L} - I]\) are all less than zero, the disease will not spread throughout the entire network system. This insight is crucial for understanding the stability and containment of epidemic spread in network structures. To ensure the epidemic threshold condition is met, the maximum eigenvalue of $\widetilde{L}$ must be 1. Given that the overall network is connected, $\widetilde{L}$ is an irreducible, non-negative square matrix. The largest absolute value of any eigenvalue of a matrix is called its spectral radius. According to the Perron-Frobenius theorem \citep{horn2012matrix}, for an irreducible, non-negative square matrix $\widetilde{L}$, the spectral radius is a real, positive, and simple eigenvalue. Furthermore, the corresponding Perron vector is strictly positive. Let $V_1$ and $V_2$ denote the Perron vectors when the spectral radius equals 1.

So, for a spectral radius equal to 1, we can write as follows:
\begin{equation}    
\begin{bmatrix}
\tau_{11}A_{11} & \tau_{12}A_{12} \\
\tau_{21}A_{21} & \tau_{22}A_{22}
\end{bmatrix}
\begin{bmatrix}
V_1 \\
V_2
\end{bmatrix}
=
\begin{bmatrix}
V_1 \\
V_2
\end{bmatrix}
\end{equation}

Based on equation no. (8) above, we can write the equations below (9 and 10).
\begin{equation}    
\tau_{11}A_{11} V_1 + \tau_{12}A_{12} V_2 = V_1
\end{equation}
\begin{equation}
\tau_{21}A_{21} V_1 + \tau_{22}A_{22} V_2 = V_2 \\
\end{equation}
Finding the value of $V_2$ from equation (10), we get
\begin{equation}
V_2 = \tau_{21}(I - \tau_{22}A_{22})^{-1}A_{21}V_1
\end{equation}
As $V_1$ is positive, it follows that $V_2 = \tau_{21}(I - \tau_{22}A_{22})^{-1}A_{21}V_1$ exists and is non-negative. By substituting the expression for $V_2$ into equation (11) and then into equation (9), we can rewrite equation (9) as $HV_1 = V_1$, where $H$ is given by:
\begin{equation}
H=\tau_{11}A_{11}+\tau_{21}\tau_{12}A_{12}(I - \tau_{22}A_{22})^{-1}A_{21}.
\end{equation}
We know that the infection process is the result of interactions between pairs of agents. It depends on the rate at which the disease is transmitted and the ability of the susceptible to receive the disease. So $\beta_{11},\beta_{12},\beta_{21},\beta_{22}$ are not completely independent on each other. Following the analysis of the infection rates in any interconnected network \citep{sahneh2013effect}, we can write.

\begin{equation}
    \beta_{12}\beta_{21}=\alpha^2\beta_{11}\beta_{22}
\end{equation}

As the recovery rates for all the agents are same, we can write the equation (13) as 
\begin{equation}
    \tau_{12}\tau_{21}=\alpha^2\tau_{11}\tau_{22}
\end{equation}
Here $\alpha$ is a positive scalar accounting for the heterogeneity of interconnection and intraconnection.
So, by substituting equation (14) in equation (12), we obtain

\begin{equation}
H=\tau_{11}{A_{11}+\alpha^2\tau_{22}\tau_{11}A_{12}(I - \tau_{22}A_{22})^{-1}A_{21}}.
\end{equation}

We can update equation (15) as $H=\tau_{11}H_{T}$ , where $H_{T}$ can be defined as
\begin{equation}
H_{T}=[A_{11} + \alpha^2\tau_{22}A_{12}(I - \tau_{22}A_{22})^{-1}A^T_{12}]
\end{equation}

As we have stated, equation (9) can be expressed in the form of $HV_1=V_1$, so from equations (15) and (16), we can  also update the expression $HV_1=V_1$ as 
\begin{equation}
\tau_{11}H_{T}V_1=V_1
\end{equation}

We said $V_1$ is a Perron vector, which is positive, and from the analysis of Sahneh et al. \citep{sahneh2013effect}, we can conclude that $H_{T}$ is an irreducible matrix. Therefore, from the Perron-Frobenius Theorem for irreducible matrices, we can conclude that  $\tau_{11}$ is the inverse of the spectral radius of $H_{T}$.
 
If we call the epidemic threshold $\tau_{11c}$, we can say that $\tau_{11c}$ is the spectral radius of the matrix $H_{T}$, which is denoted in equation (16). So, for any given infection strength $\tau_{22}$ and coupling level, if the $\tau_{11}$ value is lower than $\tau_{11c}$, we can say the epidemic will die out. 

\subsection{The threshold curve}\label{subsubsec1}

\begin{figure}[h]
  \centering
  \includegraphics[width=0.8\textwidth]{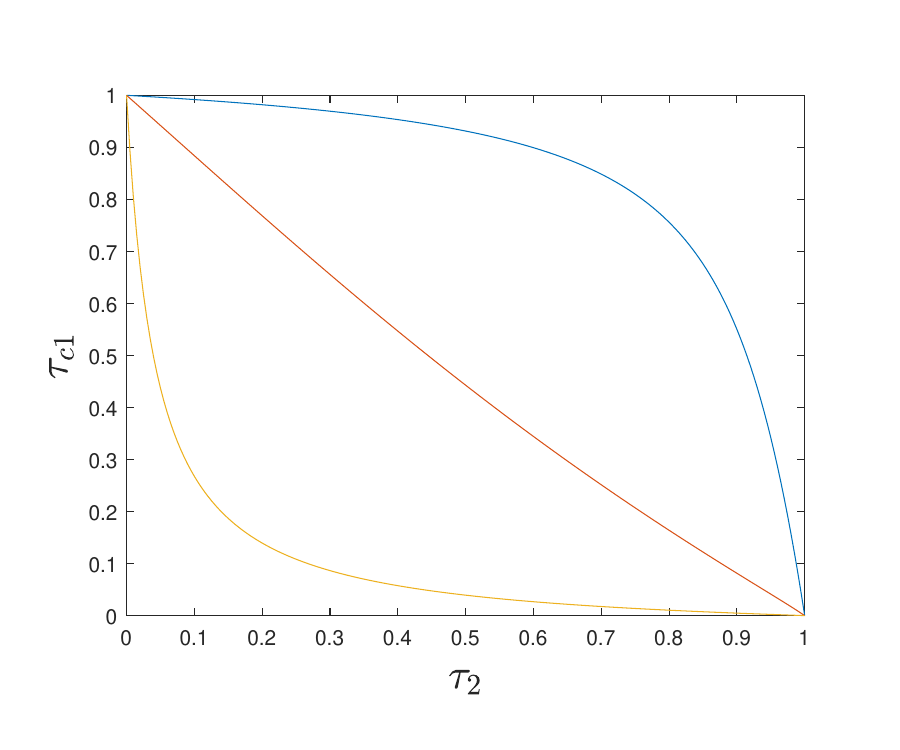}
  \caption{Graph of normalized epidemic threshold $\tau_{c1}$=$\tau_{11,c}$*$\lambda$($A_{11}$) as a function of normalized infection strength $\tau_{2}$=$\tau_{22}$*$\lambda$($A_{22}$) with different interconnection levels and $\alpha$=1. The yellow curve is for a strong level of interconnection ($\omega$=0.2), the red curve is for a medium level of interconnection ($\omega$=0.042), and the blue curve is for a low level of interconnection ($\omega$=0.01). The curve shown in this figure is for the above-mentioned small-world network topology. }
  \label{fig:2}
\end{figure}

To give an example of application of the results obtained in the previous section, we have generated two different realizations of the Watts-Strogatz (WS) model \citep{watts1998collective}. In layer 1, we have a WS network with 500 nodes, a mean node degree of 20, and a rewiring probability of 0.2. In layer 2, we have a WS network with 100 nodes, a mean degree of 4, and a rewiring probability of 0.1. We have activated all potential edges between the two layers with some probability $\omega$ for the interconnection of these two graphs. For high, medium, and low interconnection levels, the values of $\omega$ are respectively 0.2, 0.042, and 0.01. We have generated a plot of the normalized epidemic threshold $\tau_{c1}=\tau_{11,c}\lambda(A_{11})$ as a function of the normalized infection strength $\tau_{2}=\tau_{22}\lambda(A_{22})$ for different interconnection levels. The plot of epidemic threshold curves is shown in Figure \ref{fig:2}.

\begin{figure*}[h]
   \vspace*{5cm}
    \hspace*{-1cm}
       \centering
   \begin{subfigure}[b]{0.43\textwidth}
  \includegraphics[width=1.36\textwidth]{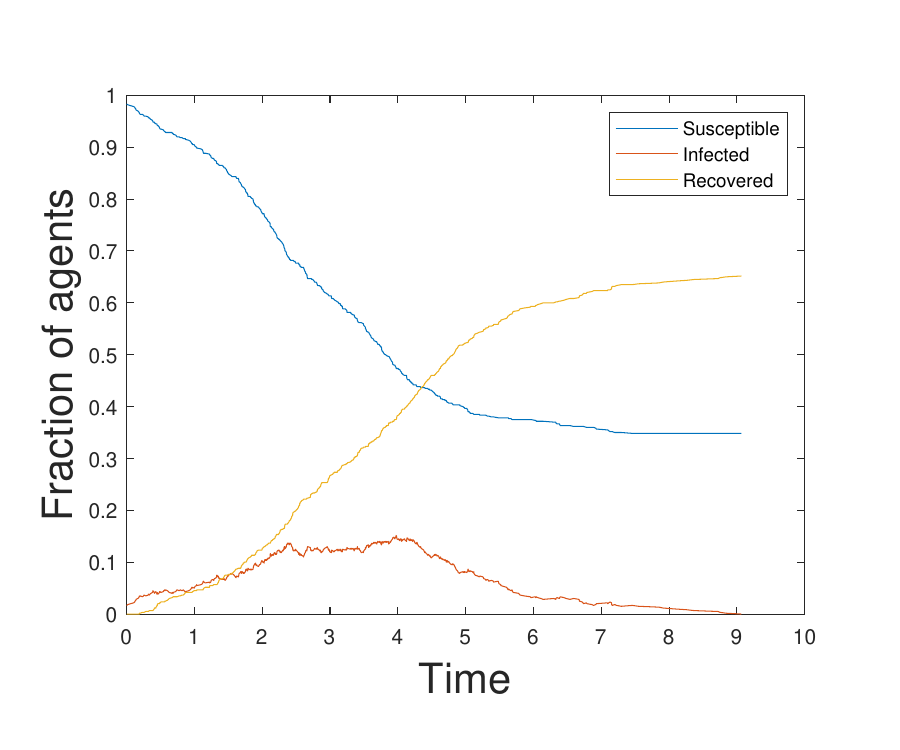}
  \label{fig:3}
   \end{subfigure}
    \hspace{1.6cm}
    \begin{subfigure}[b]{0.43\textwidth}
   \includegraphics[width=1.36
   \textwidth]{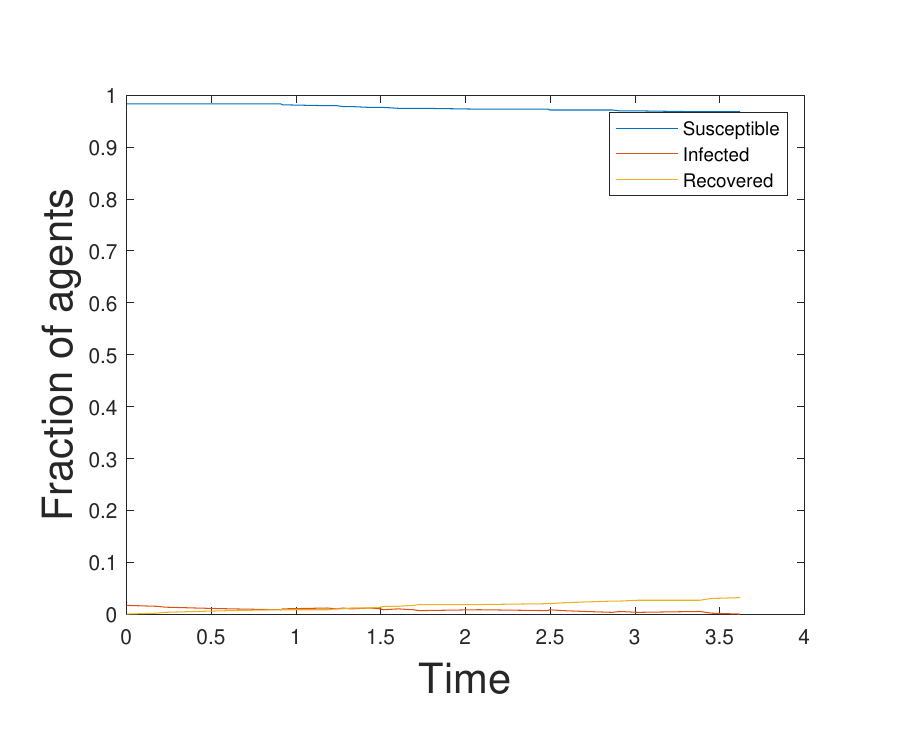}
  \label{fig:4}
  \end{subfigure}%
  \caption{Figures depicting spreading scenarios for weak and strong-level interconnections. The figure shows that when the normalized epidemic strengths are above the threshold, the infection spreads in the contact network (left). When the normalized epidemic strengths are below the threshold, the infection spreads in the contact network (right). }
    \label{fig:main}
\end{figure*}

We have taken various normalized epidemic strengths ($\tau_{c1}$ and $\tau_{2}$) and found that when the epidemic strengths are below the threshold curve, the infection is not spreading, whereas when the infection strengths are above the threshold curve, the infection is spreading. This is true for any level of interconnection. Through the following figures, we have shown how the number of agents in different compartments varies with time. The network topologies are the same Watts-Strogatz models as mentioned earlier.
In the given scenario where the normalized infection strengths ($\tau_{c1}$ and $\tau_{2}$)  for both network one and network 2 are 0.4, and the inter-network infection strengths are equal with a value of alpha being 1, as well as a recovery rate of 1 for the agents of both networks, we observe that the epidemic is spreading when the level of interconnection between the networks are strong (Figure \ref{fig:main}). Conversely, when the level of interconnection is weak, the epidemic does not spread (Figure \ref{fig:main}).

\newpage

\begin{figure}[H]
  \centering
  \includegraphics[width=0.8\textwidth]{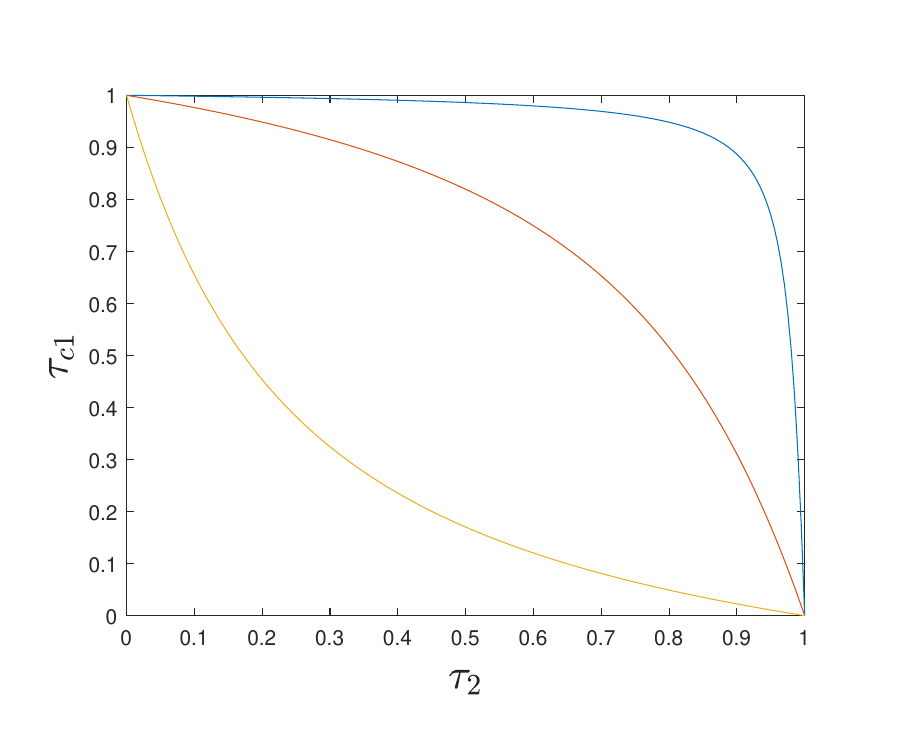}
  \caption{Graph of normalized epidemic threshold $\tau_{c1}$=$\tau_{11,c}$*$\lambda$(A11) as a function of normalized infection strength $\tau_{2}$=$\tau_{22}$*$\lambda$(A22) with different interconnection level and $\alpha$=1 for Gilbert networks. The yellow curve is for the high number of interconnections ($\omega$=0.2), the red curve is for a medium number of interconnections ($\omega$=0.042), and the blue curve is for a low number of interconnections ($\omega$=0.01). The network topologies for both networks are Erdos-Renyi.}
  \label{fig:5a}
\end{figure}

\newpage
We have also plotted the epidemic threshold curves for the Erdos-Renyi network. \citep{erdHos1960evolution} (Figure \ref{fig:5a}). The first layer is a Gilbert model with 500 nodes and a probability of interconnection of 0.02. The second layer is a Gilbert model with 100 nodes and a probability of interconnection of 0.1. For this configuration, the epidemic threshold curve is shown in Figure \ref{fig:5a}    

\begin{figure}[H]
  \centering
  \includegraphics[width=0.8\textwidth]{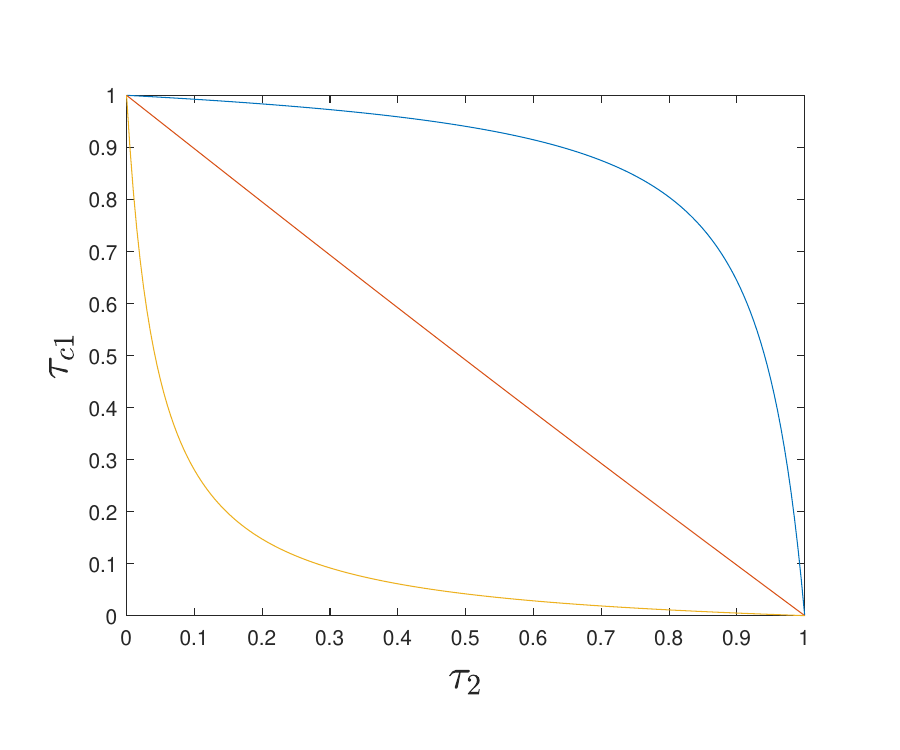}
  \caption{Graph of normalized epidemic threshold $\tau_{c1}$=$\tau_{11,c}$*$\lambda$(A11) as a function of normalized infection strength $\tau_{2}$=$\tau_{22}$*$\lambda$(A22) with different interconnection level and $\alpha$=1 for Gilbert networks. The yellow curve is for the high number of interconnections ($\omega$=0.2), the red curve is for a medium number of interconnections ($\omega$=0.042), and the blue curve is for a low number of interconnections ($\omega$=0.01). The network topologies for both networks are Erdos-Renyi.}
  \label{fig:5b}
\end{figure}

\newpage
We have also plotted the epidemic threshold curve for an interconnected network system where both the network is an Erdos-Reyni network of 500 nodes, and the probability of interconnection is 0.02. The curve is shown in Figure \ref{fig:5b} 

\clearpage

\section{Modelling Spillover}\label{sec4}

Spillover refers to transferring a virus from its usual circulating species, known as the reservoir, to a different species, known as the novel host. In the host species, the virus can either perish or undergo adaptations, potentially leading to the emergence of an epidemic. Zoonotic infectious diseases spread from animals to humans. Sixty percent of human infectious diseases are zoonotic, and seventy-five percent are emerging zoonoses \citep{salyer2017prioritizing}. Many emerging zoonotic diseases are caused by viruses, including avian influenza, rabies, and Ebola\citep{sahu2021emergence,nandi2021probability}.

\begin{figure}[H]

  \centering
  \includegraphics[width=0.8\textwidth]{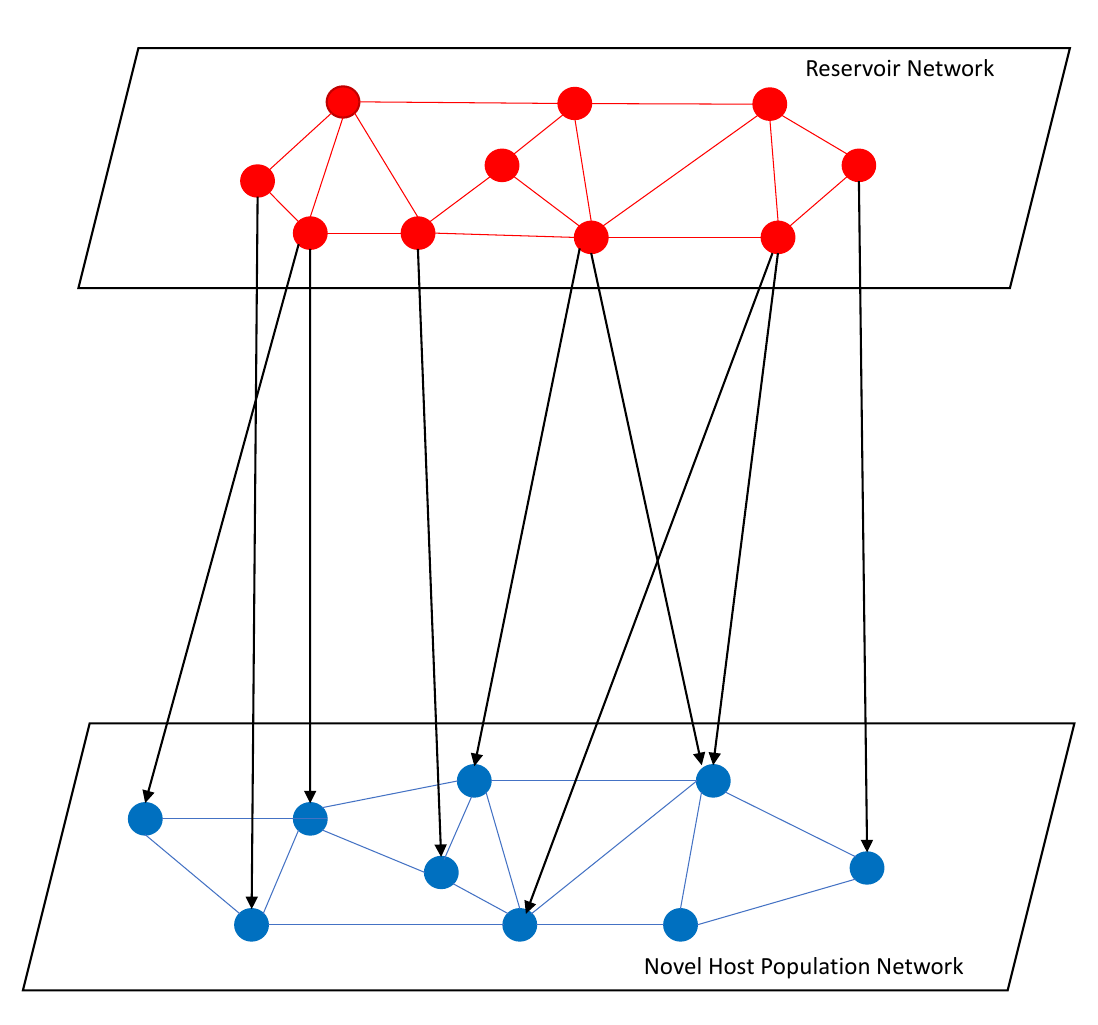}
  \caption{Spillover in networks from reservoir population to the novel host population.}
  \label{fig:6}
\end{figure}

The graphical representation in Figure \ref{fig:6} illustrates the transfer of infection, referred to as spillover. The upper network represents a reservoir network, where the nodes are highlighted in red. The lower network represents the population of a novel host, with its nodes marked in blue. The red and blue links between the agents of reservoir and novel host population network denotes the intra network connection between the agents. The black-directed links indicate the transmission of the disease from the infected agents in the reservoir network to the agents in the novel host population network.  


In a contact network where there are already some infected individuals, the spread of a disease can occur, and various factors influence the extent of the spread. However, if we keep all other properties of the network constant, it has been observed that when the initial number of infected individuals is low, the threshold for an epidemic to occur becomes higher. This means that if the initial number of infected individuals is very low, the spread of the disease can be controlled \citep{machado2022effect}.

In our work, when the number of infected individuals in the new host network exceeds a significant threshold (considered as three), it indicates an effective spillover event. Therefore, if the number of infected agents resulting from spillover in the host network is less than three, the spillover size is considered zero. Similarly, if the probability of spillover in the host network (calculated as the ratio of spillover occurrences to the total number of realizations) is less than 0.1, it is classified as a non-spillover scenario. These assumptions are in line with previous research on the stochastic epidemic model \citep{andersson2012stochastic}. 

To simulate the spillover events, we have used the stochastic simulator GEMFsim based on the  Generalized Epidemic Mean Field (GEMF) framework developed by the Network Science and Engineering (NetSE) group at Kansas State University \citep{sahneh2013generalized, sahneh2017gemfsim}. We have used the Matlab and Python version for simulation.

\subsection{Simulation results}\label{subsubsec1}

Our simulation model pertains to diseases like avian influenza, West Nile virus, and Ebola. 
In the first portion of our simulation based  study we have modeled avian influenza and West Nile virus diseases. For West Nile virus and avian influenza humans serve as dead-end hosts. Transmission occurs from birds to other birds and occasionally from birds to humans. In our simulation-based study of avian influenza and West Nile virus, we can envision a scenario in which the agents within the reservoir network represent birds while the agents in the novel host population network represent humans. Consequently, within our spillover model, we establish an effective infection rate of zero between the agents in the novel host population network. Conversely, both the effective infection rate within the reservoir network and the effective infection rate between the novel host population network and the reservoir network are non-zero. 

 The examination of data from the 2019 European study \cite{european2019annual} revealed that the percentage of infected birds can vary significantly due to various factors. With this study as a reference, we aim to establish an effective infection rate within our reservoir network such that the number of infected agents in the reservoir network falls within the range of 4 to 6 percent of the total population. Notably, the study \cite{european2019annual} found that, on average, 5.6 percent of birds were infected with Highly Pathogenic Avian Influenza (HPAI) in eight European countries.

In this specific situation, both the novel host network and reservoir network are Erdos-Renyi networks consisting of 1000 nodes and 3255 links. Throughout the simulation, the infection rates have been set as follows. The infection rate within reservoir network $\beta_{22}$ is set to 0.15 in order to obtain the expected number of infected agents to be within 4 to 6 percent of the total population as reported in \cite{european2019annual}. 

Now, the probability that a human is infected by a bird can vary a lot based on different factors like virus strain, immunity of people, infection-spreading capability of birds, etc. However, based on a study conducted on Egyptian people, the probability of a person being infected by a poultry bird is 0.02 \citep{gomaa2015avian}. Another study conducted in European countries claims that the risk of avian influenza for general people is low, and for occupationally exposed people, the risk of infection is low to medium \citep{european2023avian}. 
Based on these studies, the infection rate between agents of the reservoir network and the novel host population network  $\beta_{12}$ is set equal to 0.02 such that the probability of infection for humans from birds is 0.02, and the rate is also low to moderate. Finally, we set the infection rate within the novel host population network $\beta_{11}$ equal to 0, assuming human-to-human disease transmission is very limited. The recovery rate $\delta$ for all agents is assumed to be equal to 1 to maintain a proper effective infection rate. The fraction of the total number of interpopulation links between the two networks has ranged from 0.0001 to 0.03. To establish an interpopulation link, we randomly select one node from the reservoir network and one node from the novel host population network, creating a directed link from the reservoir network to the novel host population network.

\begin{figure}[H]
  \centering
 \includegraphics[width=1.1\textwidth]{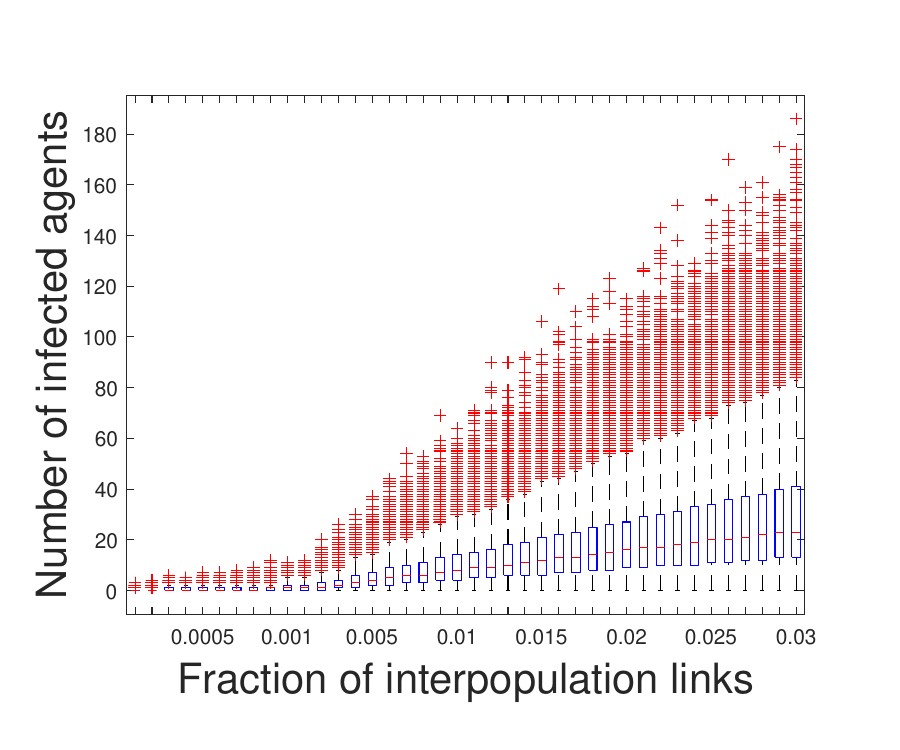}
  \caption{The number of infected agents in novel host population network due to spillover with the fraction of total possible links between novel host population network and reservoir network. Every 16000 realizations have been shown in the boxplot.}
  \label{fig:7}
\end{figure}

To measure the spillover size, 16000 simulations have been conducted for each fraction of all possible links. For each specific fraction of links, a box plot has been generated to depict the spillover size (Figure \ref{fig:7}). Ten random nodes were chosen as initially infected nodes in the reservoir network. This random selection process was repeated for each realization, resulting in different sets of initially infected nodes for each iteration. 

\begin{figure}[H]
  \centering
  \includegraphics[width=0.8\textwidth]{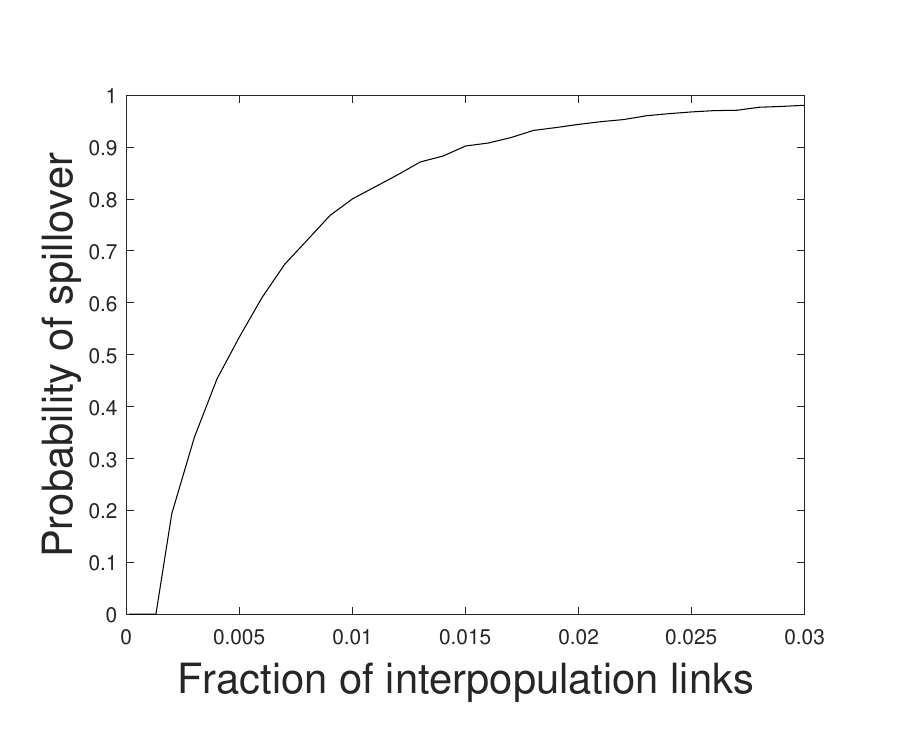}
  \caption{Probability of spillover in novel host population network with inter-population link density. All other parameters are kept constant. }
  \label{fig:8}
\end{figure}

Based on the assumptions about spillover mentioned earlier in this paper, we have generated a plot between the probability of spillover in the host network and the fraction of links between the reservoir and the host network (Figure \ref{fig:8}).

It is important to mention that when the fraction of links is smaller than 0.0013, the probability of spillover remains very low. However, beyond this threshold, there is a noticeable increase in the probability of spillover. Therefore, we can conclude that the fraction of 0.0013 of the total possible inter-network links acts as a clear threshold, indicating a distinct phase transition.

In the next scenario, the number of links between the novel host population network and reservoir network has been kept constant at 1000. However, the inter-network infection rate ($\beta_{12}$) has been systematically reduced from 0.7 to 0.005. The infection rate within the novel host population network  ($\beta_{11}$) is set to zero, while the infection rate within the reservoir network ($\beta_{22}$) remains at 0.15. The novel host network and the reservoir network maintain the same topology as previously mentioned. It is assumed that ten random nodes are infected during each iteration.

\begin{figure}[H]
  \centering
  \hspace*{-1.1cm}
  \includegraphics[width=0.8\textwidth]{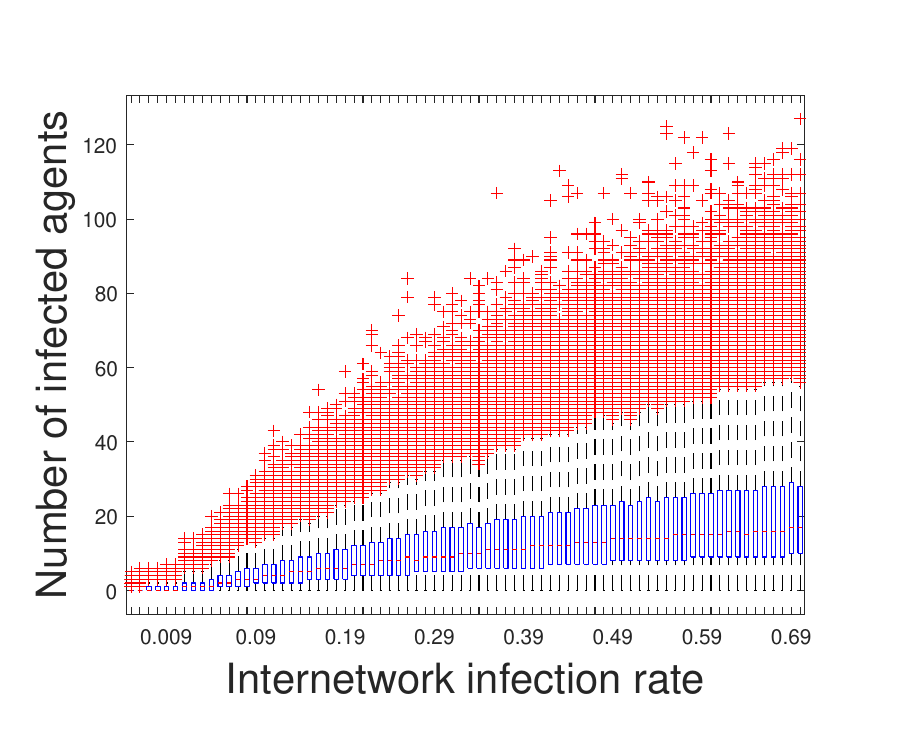}
  \caption{The number of infected agents in novel host population network with the inter-network infection rate between novel host population network and reservoir network. Every 16000 realizations have been shown in the boxplot. These boxplots are generated in Matlab}
  \label{fig:9}
\end{figure}

To analyze the impact of different inter-network infection rates on spillover size, 16000 realizations have been recorded for each value of inter-network infection rate. For each inter-network infection strength, a box plot has been generated to represent the spillover size visually (Figure \ref{fig:9}). 

\begin{figure}[H]
  \centering
  \includegraphics[width=0.8\textwidth]{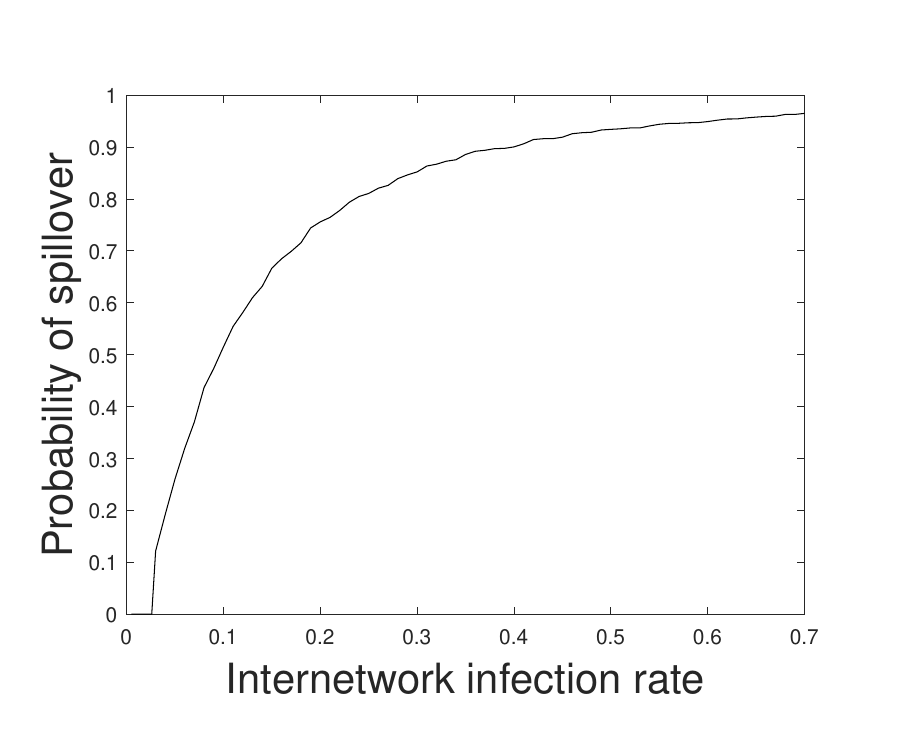}
  \caption{Probability of spillover with inter-network infection rate keeping all other parameters constant.}
  \label{fig:10}
\end{figure}

Based on the previously mentioned assumptions about spillover, we draw a curve between the probability of spillover and inter-network infection strength between two networks (Figure \ref{fig:10}).

If we observe Figure \ref{fig:10} closely, we can see that the inter-network infection rate close to 0.026 acts as a threshold between a low and significantly high probability of spillover. A clear phase transition in spillover probability also exists for the change in inter-network infection rate.

\begin{figure}[H]
  \centering
  \includegraphics[width=0.8\textwidth]{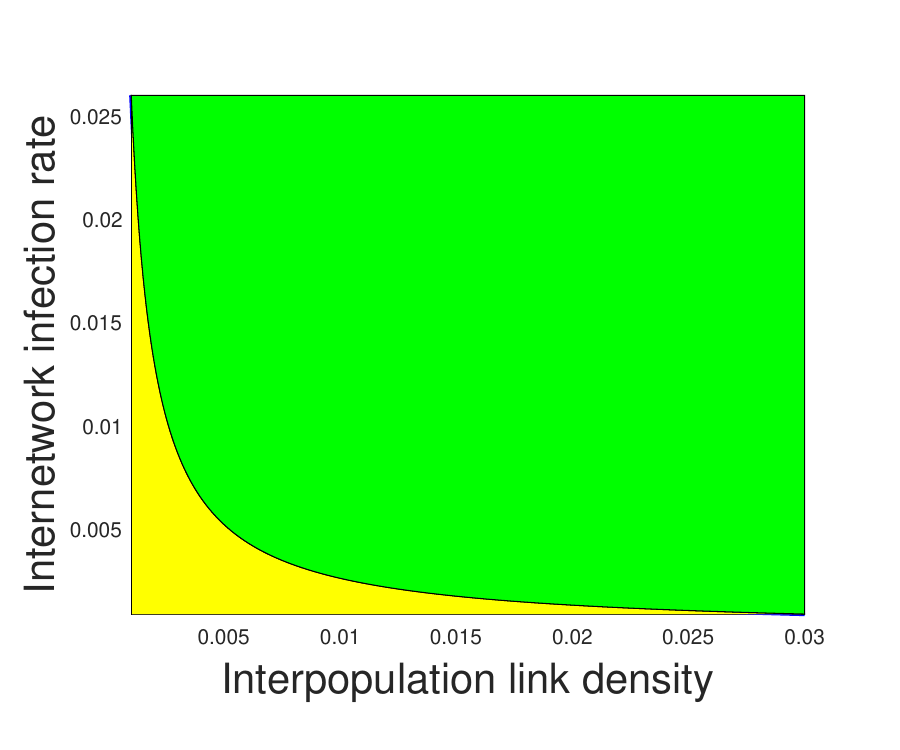}
  \caption{ Spillover region separation based on inter-population links and inter-network infection strength. The green-colored zone is the significant spillover zone, while the yellow zone is the non-significant spillover zone.}
  \label{fig:11}
\end{figure}

When all other factors remain unchanged, and the fraction of total possible links between different populations exceeds approximately 0.001 in a fixed network structure, we can observe a relationship between the strength of infection between networks ($\beta_{12}$) and the total fraction of inter-population links between populations. This relationship confirms the validity of our assumption regarding the extent of spillover events. In fact, the strength of infection between networks ($\beta_{12}$) and the proportion of links between populations are inversely proportional to each other. Based on this correlation, we can represent a rectangular hyperbola curve where either of these two quantities can be plotted on the x-axis and the other on the y-axis (Figure \ref{fig:11}). The region above this curve (shown in green in Figure \ref{fig:11}) corresponds to a significant spillover area, while the region below the curve represents a minor spillover area (shown in yellow in Figure \ref{fig:11}).

Furthermore, we conducted the same experiment with scale-free networks \citep{barabasi2003scale} in both layers. There were no significant changes in the results compared to the previous findings when the inter-network links between the host and reservoir networks were randomly connected. 

\begin{figure}[H]
  \centering
  \includegraphics[width=0.8
  \textwidth]{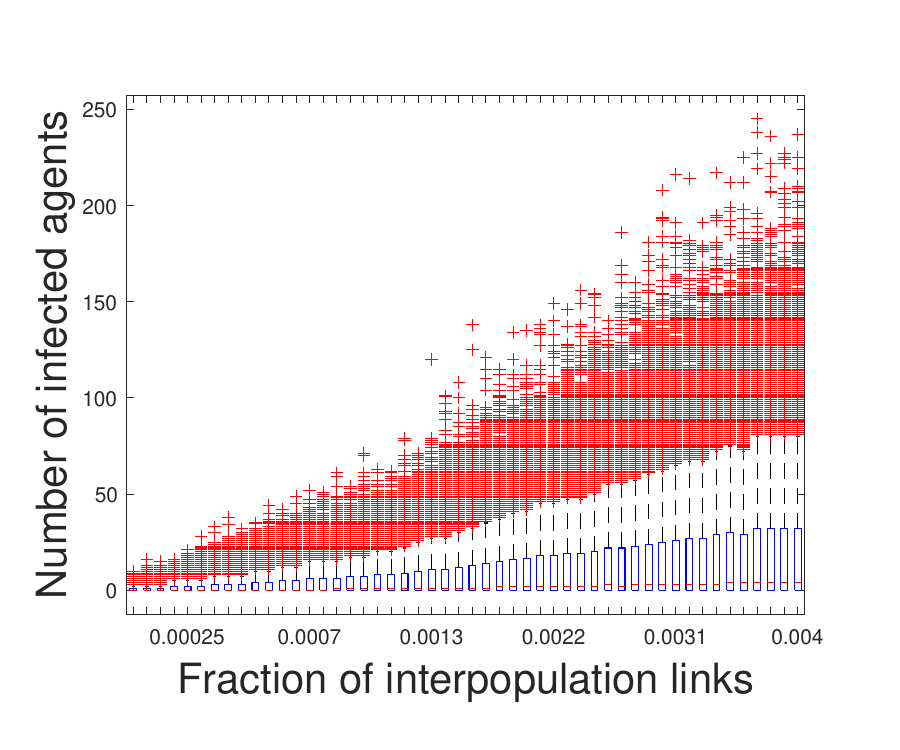}
  \caption{The inter-population links are connected with only the hubs of the reservoir network. The spillover size for each fraction of links has been shown in boxplot}
  \label{fig:12}
\end{figure}

\begin{figure}[H]
  \centering
  \includegraphics[width=0.8\textwidth]{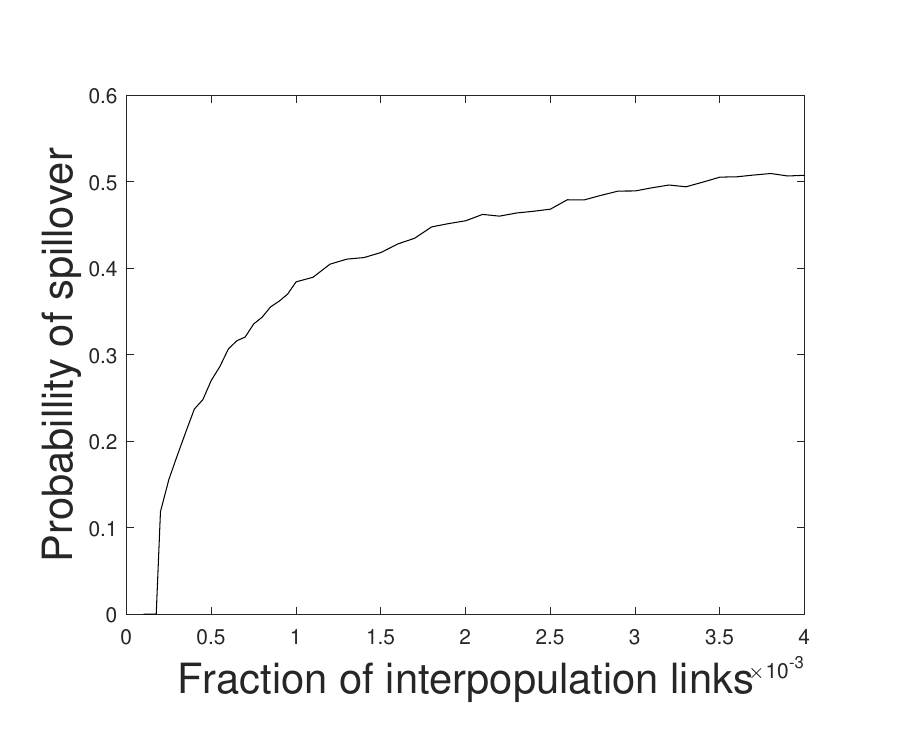}
  \caption{The phase transition curve when the inter-population links between novel host population network and reservoir network are connected with the hubs of reservoir network. The phase transition occurs for a very low fraction of links. }
  \label{fig:13}
\end{figure}

However, when the agents of the novel host network were connected with the hubs of the reservoir network using the same inter-network infection rate, a phase transition occurred at a significantly lower fraction of links compared to the previous results shown in Figure \ref{fig:8}. In this case, we had a scale-free network consisting of 1000 nodes and 2957 links for both the novel host population network and the reservoir network. The infection rate among the agents of the reservoir network was 0.1, and the inter-network infection rate was 0.02. We distributed the number of links equally among five hubs. However, the fraction of links varied from 0.0001 to 0.0039, and for each fraction of links, we conducted 16000 realizations to illustrate the size of spillover, shown in the boxplot (Figure \ref{fig:12}). We observed the same phase transition in the spillover probability as before, but it occurred when the fraction of inter-population links was much smaller, approximately 0.00018 (Figure \ref{fig:13}).

In our study, we expanded our simulations to include small-world networks (based on Watts and Strogatz's model) and observed that the same phase transition phenomenon also exists for these networks. 

To summarize, the primary focus of our simulations was on exploring the spillover threshold. To achieve this, we selected random networks, scale-free networks, and several small-world networks with varying rewiring probabilities (ranging from 0.01 to 1). All of these networks contained the same number of nodes (1000), while the small-world networks had a slightly higher number of links.

During the course of the experiments, we carefully controlled the infection rate to ensure that the average number of total infected individuals within the reservoir network remained consistently stable. Across all mentioned network topologies, this average ranged around 5 percent throughout the entire infection process. Additionally, the inter-network infection strength ($\beta_{12}$) was kept constant. Our goal was to measure the number of inter-population links required to reach the spillover threshold. For this threshold, we validated our previous assumption, stating that approximately 0.1 fractions of all realizations would exhibit a spillover non-negligible spillover size.

Our findings revealed that the number of links needed to reach the spillover threshold was the lowest for scale-free networks and the highest for regular lattice networks. As for the small-world networks, an increase in the rewiring probability resulted in a decrease in the number of inter-population links required to reach the spillover threshold. This behavior can be attributed to the emergence of the small-world effect, which promotes the spread of epidemics \citep{liu2015epidemics}.

In the case of scale-free networks, the presence of hubs, which are highly connected nodes, contributed to an enhanced spread of epidemics within the network.

\begin{figure}
    \centering
  \includegraphics[width=1\textwidth]{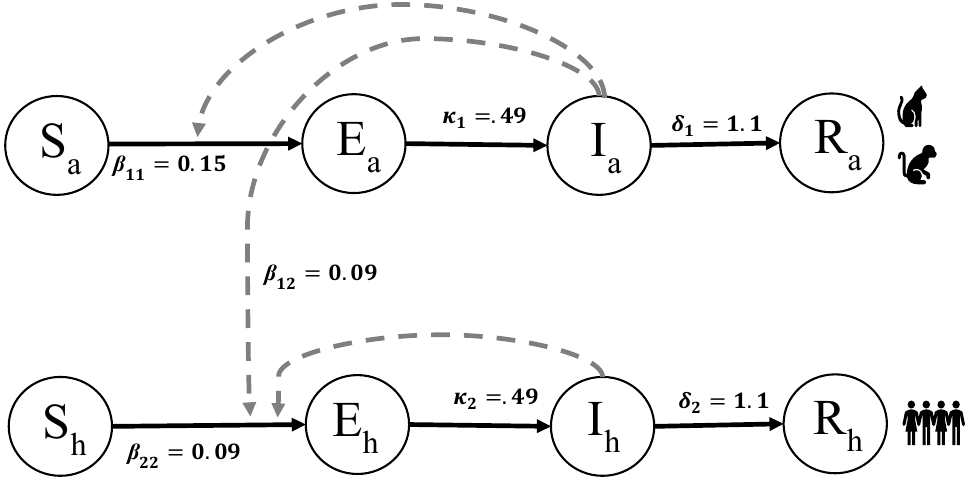}
    \caption{ SEIR  Transition Graph for Ebola Across Interconnected Networks. The diagram illustrates the transition rates between SEIR states, \textit{a} and \textit{h} subscript denote to animal and human agents, respectively. Parameters $\beta_{11}$ and $\beta_{22}$ represent the infection rates between the agents of the novel host population network and reservoir network, respectively. and $\beta_{12}$ denotes infection rate between the agents of the reservoir network and novel host population network ; $\kappa_1$ and $\kappa_2$ indicate the incubation rates of animals and humans respectively that lead to infectious states; $\delta_1$ and $\delta_2$ are the recovery rates of animals and humans respectively. Dashed arrows show the inducing state for transitions.
}
    \label{fig:14}
\end{figure}

In another scenario, our simulation-based study has been expanded to discuss one of the most notorious zoonotic diseases, Ebola, which is prevalent in West African countries. Unlike other infections where humans serve as dead-end hosts, Ebola can also spread through human networks. The Susceptible-Exposed-Infected-Recovered (SEIR) model is particularly apt for modeling Ebola due to its ability to incorporate the incubation period before the disease becomes transmissible. This model enhances the SIR model by introducing an 'Exposed' state that accounts for the incubation period. The transition between states in the SEIR model, as shown in \ref{fig:14}, includes the incubation rate ($\kappa$), which is incubation rate related the time it takes for an exposed individual to become infectious.

\begin{figure}[t]
  \centering
  \includegraphics[width=0.8\textwidth]{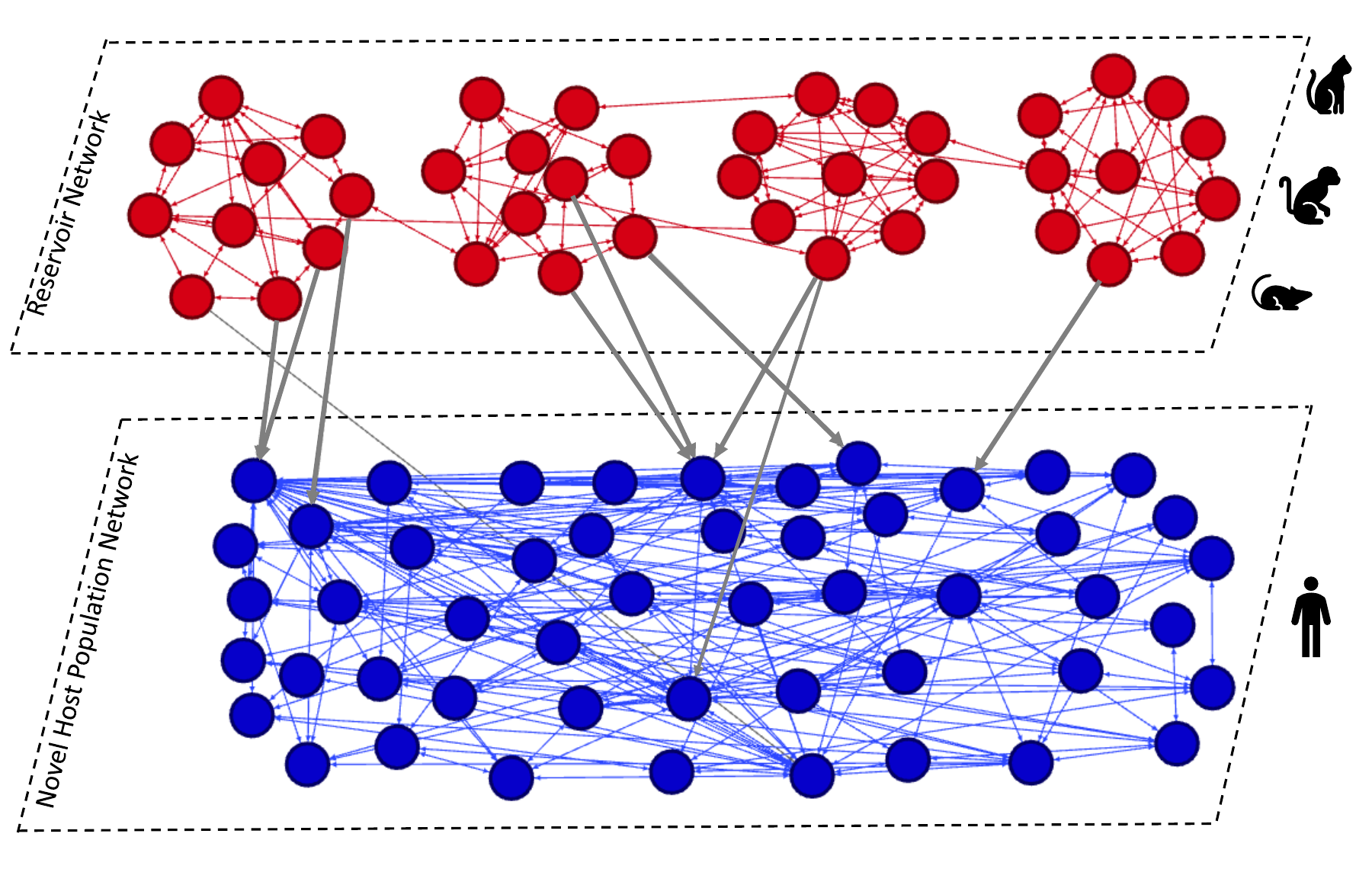}
  \caption{Illustration of an example for two-layer network structure used in the simulation: (a) represents the contact network among animals, (b) represents the contact network among humans, and the grey arrows denote the directed transmission links from animals to humans.}
  \label{fig:15}
\end{figure}

In this analysis, we examine a hypothetical yet plausible scenario for the transmission of the Ebola virus from animal reservoirs to human populations. To enhance the realism of our model, we constructed two distinct network layers. The first layer, representing the animal population, is modeled as a stochastic block network comprising 2000 nodes and 5079 undirected edges. This network simulates 20 distinct animal herds, where each block represents a particular herd, each connected to at least two other herds, mimicking the potential inter-herd transmission dynamics. For the human population layer, we employed a Barabási-Albert (BA) network to represent a small village setting. The BA network, renowned for its scale-free properties, accurately captures the intricate social structures present in human communities, featuring 5000 nodes and 14991 undirected edges. In Figure \ref{fig:15}, an example  of these two networks structure is presented. We extracted the parameters for the SEIR model from the cases covered in \cite{drake2015transmission} as Uganda 2000 case and adjusted the rates to accurately reflect the network structure for both humans and animals. The rates were adapted based on the average degree of nodes $\langle k \rangle$. The parameters used for this experiment are detailed in Figure \ref{fig:14}.\\

\newpage

\begin{figure}[H]
  \centering
  \includegraphics[width=1\textwidth]{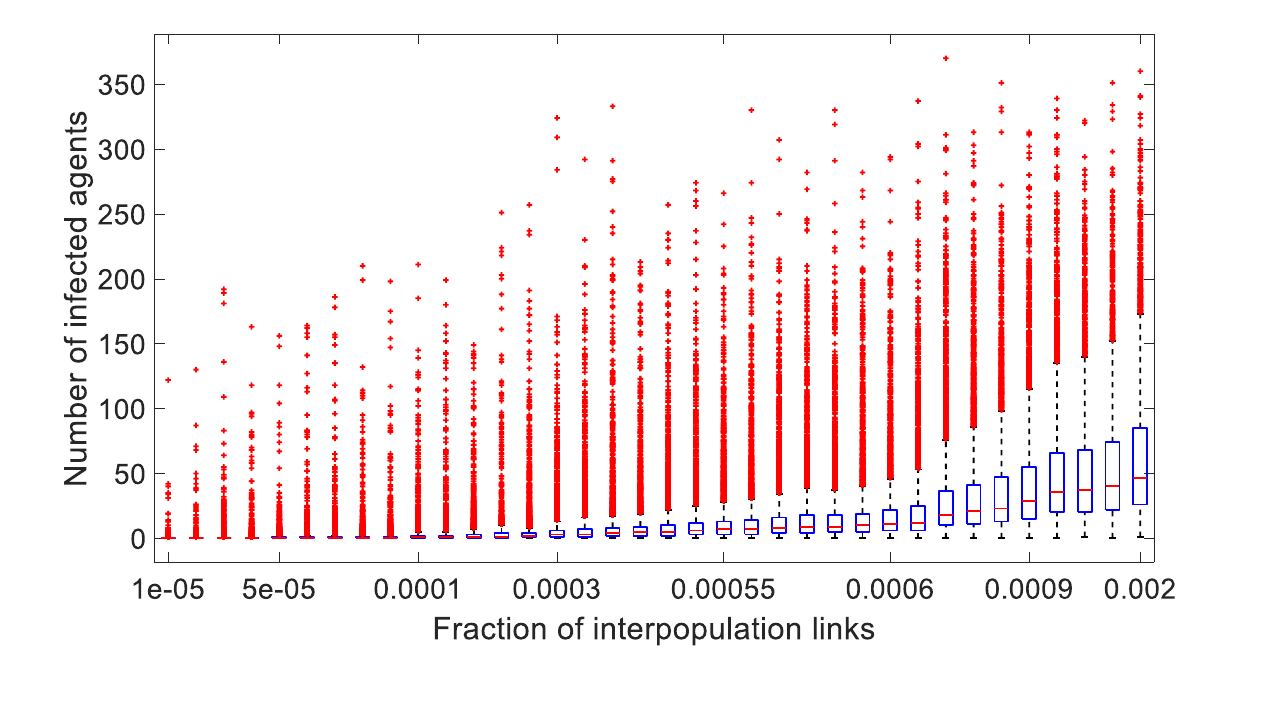}
  \caption{The boxplot of number of infected agents in novel host population network for simulation-based study of Ebola. The number of infected agents in the novel host population network increases with the fraction of the possible links between the novel host population network and the reservoir network.}
  \label{fig:16}
\end{figure}

\begin{figure}[H]
  \centering
  \includegraphics[width=1\textwidth]{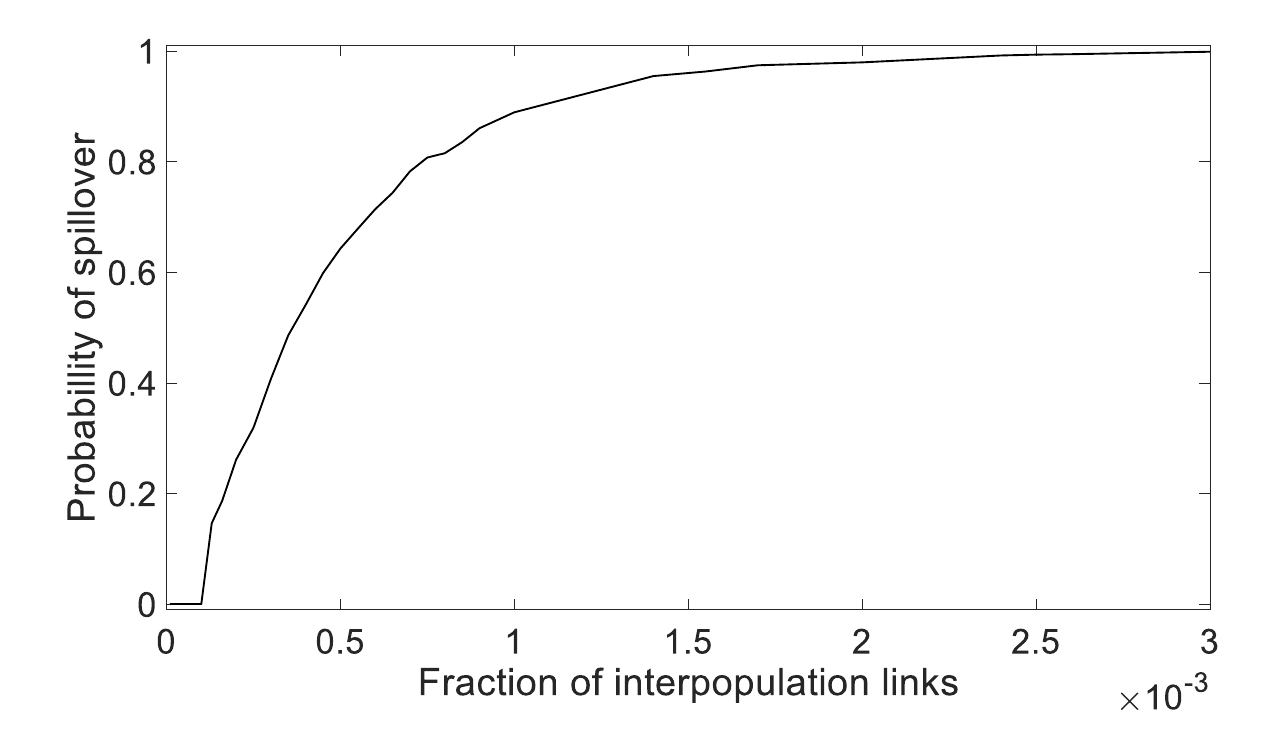}
  \caption{Probability of spillover in host population network for simulation-based study of Ebola}
  \label{fig:17}
\end{figure}

We ran 16000 simulations to examine the effects of varying the number of links between humans as hosts and animals as networks as we increase the links from 100 to 30,000 edges. Figure \ref{fig:16} shows the number of infected agents in the novel host population network for a specific fraction of interpopulation links. The data from a total of 16,000 simulation realizations is represented using a boxplot. The probability of spill-over with respect to the fraction of interpopulation links is also presented in Figure \ref{fig:17}. In this simulation-based study, we can also see a clear phase transition for the probability of spillover as the number of interpopulation links increases. A fraction of up to 0.0001 of the total possible number of links acts as a threshold for this phase transition curve.

Based on the simulation-based study of the above-mentioned zoonotic diseases, we can conclude that there will always be a clear phase transition in the probability of spillover as the number of interpopulation links increases, regardless of whether humans act as dead-end hosts.

\newpage
\section{Conclusions}\label{sec5}

Our work investigated two situations for SIR spreading dynamics in interconnected networks. In the first case, we have a situation where the initial infection could be present in any of the two interconnected networks. The infection will not spread in interconnected networks if the normalized infection strengths are below the threshold curve for any level of interconnection. However, in the second situation, the infection in the novel host population comes only from the reservoir network. For avian influenza and West Nile virus modelling, we conducted experiments by changing the number of inter-population links and by changing the inter-network infection rate, keeping all other parameters constant. Experiments have been conducted on Erdos-Renyi networks, Barabasi-Albert networks, and Watts-Strogatz networks. The inter-population links were distributed randomly. In all the cases, we obtained a phase transition in the probability of spillover. While keeping all other network parameters constant, we obtain a relationship between the inter-population link density and inter-network infection strength, based on which we get two regimes of significant and minor spillover zones. When the inter-population links are connected only with the hubs, we get a phase transition at a very low number of links while keeping all other parameters constant. It can be concluded that hubs of scale-free networks also play a crucial role in spillover. For small-world networks, we found that the number of interpopulation links for the spillover threshold decreases with the increase in rewiring probability. The spillover threshold, when considering a fixed number of infections in the reservoir population, follows a pattern of highest to lowest interpopulation links: Regular Lattice networks have the highest number of links, followed by Small-World networks, Erdos-Renyi networks, and finally Scale-Free networks with the lowest number of links. Finally, when we expanded our simulation-based study to model Ebola, where human-to-human infection spreading is significant, we observed a clear phase transition in the probability of spillover.

For future work,  this study can be enhanced by developing a theoretical model to determine the spillover threshold for both inter-population links and inter-network infection strength. Additionally, simulations using various network topologies and inter-population link patterns can provide diverse results. These simulation outcomes can be used to train a neural network for predicting spillover in new scenarios.
 \section*{Abbreviations}
 
 \begin{itemize}
     \item SIR: Susceptible-Infected-Recovered
     \item SIS: Susceptible-Infected-Susceptible
     \item SEIR: Susceptible-Exposed-Infected-Recovered
     \item WS: Watts-Strogatz
     \item GEMF: Generalized Epidemic Mean Field
     \item NetSE: Network Science and Engineering
     \item HPAI: Highly Pathogenic avian Influenza
     \item BA: Barabási-Albert  
     
 \end{itemize}

\section*{Declarations}

\begin{itemize}

\item Ethics approval and consent to participate:- Not applicable

\item Consent for publication:- Authors give their consent for publication.

\item Funding:- This research was supported by the U.S. Department of Agriculture under Grant Number 2022-67015-38059. 

\item Competing interest:- The authors confirm that they do not have any competing interests to report.

\item Data and Code Availability:- The data and code used in this study are available upon request.
 
\item Authors' contributions:-

S.D. and C.S. contributed to conceptualization and methodology, S.D. prepared the original draft, C.S. contributed to supervision,  project administration, and securing funding. S.D. and C.S. reviewed and edited the manuscript. M.H.S. contributed  to the simulation work presented in this paper and assisted in the verification of the final manuscript, ensuring the consistency of data presented."

\item Acknowledgements:- We would like to acknowledge the Network Science and Engineering (NETSE) group at Kansas State University.

\end{itemize}

\bibliography{references}
\end{document}